\RequirePackage{fix-cm}
\documentclass[smallextended]{svjour3}
\smartqed

\journalname{Optimization Letters}
\usepackage{hyperref}
\hypersetup{
    unicode = false,
    pdftoolbar = true,
    colorlinks = true,
    linkcolor = blue,
    citecolor = blue,
    filecolor = black,
    urlcolor = blue,
    breaklinks = true
}

\usepackage{graphicx}
\usepackage{amssymb}
\usepackage{epstopdf}
\usepackage{amsmath}
\usepackage{amsfonts}
\usepackage{amssymb}
\usepackage{amsthm}
\usepackage{enumerate}
\usepackage{framed}
\usepackage{hyperref}
\usepackage{geometry}
\usepackage{makeidx}
\usepackage{epic}
\usepackage{xspace}
\usepackage{hhline}
\usepackage{curves}
\usepackage{multicol}
\usepackage{multirow}

\newcommand{\Proof}{\noindent {\bf Proof:} }
\newcommand{\R}{{\rm I\!R}}

\newcommand{\dom}{\mathrm{dom}}
\newcommand{\cl}{\mathrm{cl}}

\newcommand{\argmax}{\mathrm{argmax}} 

\begin{document}

\title{Compositions of Convex Functions and Fully Linear Models}
\author{W. Hare \thanks{Mathematics, University of British Columbia, Okanagan Campus, 3187 University Way, Kelowna, BC, Canada.       {\tt warren.hare@ubc.ca}}}

\maketitle

\begin{abstract} Derivative-free optimization (DFO) is the mathematical study of the optimization algorithms that do not use derivatives.  One branch of DFO focuses on model-based DFO methods, where an approximation of the objective function is used to guide the optimization algorithm.  Proving convergence of such methods often applies an assumption that the approximations form {\em fully linear models} -- an assumption that requires the true objective function to be smooth.  However, some recent methods have loosened this assumption and instead worked with functions that are compositions of smooth functions with simple convex functions (the max-function or the $\ell_1$ norm).  In this paper, we examine the error bounds resulting from the composition of a convex lower semi-continuous function with a smooth vector-valued function when it is possible to provide fully linear models for each component of the vector-valued function.  We derive error bounds for the resulting function values and subgradient vectors.\end{abstract}



\section{Introduction}\label{introduction}

Derivative-free optimization (DFO), the mathematical study of the optimization algorithms that do not use derivatives, has become an important aspect in modern optimization research.  One of the most common uses of DFO is for the optimization of an objective function that results from a computer simulation.  As simulation-based research becomes more prevalent, the need for DFO methods is likely to continue to increase.

It is important to note that, while DFO algorithms do not use derivatives, this is not the same as saying the objective function is not differentiable.  For example, if numerical integration is being used, then the objective function is very likely to be differentiable, but derivatives may be very difficult to obtain analytically.  This insight has inspired a collection DFO methods that are based on approximating the objective function with a model function.  The gradients, or even second derivatives, of the model function can then be used to guide the optimization algorithm.  Examples of such methods include \cite{ConnScheinbergToint1997} \cite{Powell02}  \cite{CustVicente2007} \cite{Powell2008} \cite{ConnScheinbergVicente2009}  \cite{HareNutini2013} \cite{HareLucet2014} \cite{BauschkeHareMoursi2015} (among many others).

Intuitively, in order to prove convergence for a model-based DFO method, one requires that the model functions are sufficiently accurate approximations of the true objective function.  In recent research, this notion has been mathematically captured in the definition of {\em fully linear models} \cite{conn-scheinberg-vicente-2008} \cite{WildShoemaker2011}.  Fully linear models are models that approximate the true objective function (near the incumbent solution) in a manner similar to a first-order Taylor expansion (an exact definition appears in Section \ref{ssec:notation}).

To define fully linear models, we must assume that the true objective function is smooth, $f \in \mathcal{C}^1$.  Thus, a standard assumption in model-based DFO research is smoothness of the objective function.  If the objective function is smooth, and the gradient mapping is locally Lipschitz, then a number of methods have been developed for constructing fully linear models \cite{Powell2003} \cite{conn-scheinberg-vicente-2008} \cite{connscheinbergvicente2008b} \cite{Regis2015} (and references therein). 
It is quite remarkable that some research has moved away from assuming $f \in \mathcal{C}^1$.  For example, in \cite{HareNutini2013}, it is assumed that the true objective function takes the form $f = \max\{ F_i : i = 1, 2, ..., m \}$, where each $F_i \in \mathcal{C}^1$.   Clearly, in this situation $f \notin \mathcal{C}^1$.  However, if it is possible to provide fully linear models for each $F_i$, then it is still possible to create a convergent algorithm \cite{HareNutini2013}. (An application in seismic retrofitting design that fulfills this assumption is examined in \cite{BigdeliHareNutiniTesfamariam-2016}.) Another example is \cite{LarsonMenickellyWild2015}, where the objective function takes the form $f = \sum_{i=1}^m |F_i|$.  Again, while $f \notin \mathcal{C}^1$, if it is possible to provide fully linear models for each $F_i$, then a convergent algorithm can still be created.

Notice, in both \cite{HareNutini2013} and \cite{LarsonMenickellyWild2015}, the objective function is a composition of a convex lower semi-continuous (lsc) function with a smooth vector-valued function.  In this paper, we examine the error bounds for the composition of a convex lsc function with a smooth vector-valued function under the assumption that it is possible to provide {\em fully linear models} for each component of the vector-valued function.  This follows a similar vein to the recent work by Regis, \cite{Regis2015}, which explores calculus rules for the simplex gradient.  (Regis's work includes, sum rules, product rules, and quotient rules; however composition rules are unexplored.)  The results herein encompass the error bound results of \cite{HareNutini2013} and \cite{LarsonMenickellyWild2015}, and also makes way for other composition functions, such as an augmented Lagrangian penalty function, $\sum_{i=1}^m (\max\{0, F_i(x)\})^2$.

We find that, under mild assumptions, the error bound in the function value retains the same order, over the same neighbourhood, as the original error bound (Section \ref{functionvalueapproximations}).  We also find that, under mild assumptions, the error bound for the first order information (subdifferentials) also retains the same order as the original error bound, but only at the focal point of the fully linear models (Section \ref{subdifferentialapproximations}).  In addition, we provide some simple examples that demonstrate what can go wrong if the assumptions are removed.

\section{Fully Linear Models and other Notation}\label{ssec:notation}

The principal focus on this paper is the behaviour of {\em fully linear models}, formally defined as follows.

\begin{definition}[Fully Linear Models]
\label{definition:fullylinearatx}
Given $f \in \mathcal{C}^1$, $\bar{x}\in \R^n$, and $\bar{\Delta}>0$, we say that $\{\tilde f_\Delta \in  \mathcal{C}^1 : \Delta \in (0, \bar{\Delta})\}$ are {\em fully linear models of $f$ at $\bar x$} if there exists scalars $\kappa_f >0$ and $\kappa_g >0$ such that for all $\Delta \in (0, \bar{\Delta})$
	\begin{equation}\label{eq:fullylinear}\begin{array}{rcll}
		|f(y) - \tilde{f}_\Delta(y)| &\leq & \kappa_f \Delta^2 & \qquad \mbox{ for all }~y \in B_\Delta(\bar x), \\
		\mbox{ and } \quad
		\|\nabla f(y) - \nabla \tilde{f}_\Delta(y)\| &\leq&  \kappa_g \Delta & \qquad \mbox{ for all }~y \in B_\Delta(\bar x).
	\end{array}\end{equation}
We say $\{\tilde f_\Delta \in  \mathcal{C}^1 : \Delta \in (0, \bar{\Delta})\}$ are fully linear models of $f$ at $\bar x$ {\em with constants  $\kappa_f >0$ and $\kappa_g >0$}, if conditions \eqref{eq:fullylinear} are satisfied for these constants.
\end{definition}

Recall a function is convex if its epi-graph ($\mathrm{epi}(f) = \{(x,\alpha) : \alpha \geq f(x)\}$) is a convex set \cite[Def 2.1, Prop 2.4]{rockafellar-wets-1998}, and a function is lower semi-continuous (lsc) if its epi-graph is a closed set \cite[Def 1.5, Thm 1.6]{rockafellar-wets-1998}.  As we are concerned with the composition of a convex lsc and a smooth vector-valued function $F : \R^n \rightarrow \R^m$, we must also define fully linear models for vector-valued functions.

\begin{definition}[Fully Linear Models for Vector Valued Functions]
\label{definition:fullylinearvector}
Given $\bar{x} \in \R^n$, $\bar{\Delta}>0$, and
    $$F :\begin{array}{lcl}
    \R^n &\rightarrow& \R^m \\
    F(x) &\mapsto& (F_1(x), F_2(x), ..., F_m(x)),
    \end{array}$$
such that $F_i : \R^n \rightarrow \R$ is $\mathcal{C}^1$ for each $i = 1, 2, ..., m$, we say that $\{\tilde F_\Delta : \Delta \in (0, \bar{\Delta})\}$, with
    $$\tilde F_\Delta :\begin{array}{lcl}
    \R^n &\rightarrow& \R^m \\
    \tilde F_\Delta(x) &\mapsto& (\tilde{F}_{1,\Delta}(x), \tilde{F}_{2,\Delta}(x), ... \tilde{F}_{m,\Delta}(x)),
    \end{array}$$
are {\em fully linear models of $F$ at $\bar{x}$} ({\em with constants $\kappa_{F} > 0$ and $\kappa_{G} > 0$}) if $\{\tilde F_{i,\Delta} : \Delta \in (0, \bar{\Delta})\}$ are fully linear models of $F_i$ at $\bar{x}$ (with constants $\kappa_{F} > 0$ and $\kappa_{G} > 0$) for each $i=1, 2, ..., m$.
\end{definition}

Henceforth, we will use the notation $F(x) = (F_1(x), F_2(x), ..., F_m(x))$ and \\ $\tilde F_\Delta(x) = (\tilde{F}_{1,\Delta}(x), \tilde{F}_{2,\Delta}(x), ... \tilde{F}_{m,\Delta}(x))$.

For a vector-valued function, fully linear implies that for each $i=1, 2, ..., m$, there exists scalars $\kappa_{F_i} > 0$ and $\kappa_{G_i} > 0$, satisfying conditions \eqref{eq:fullylinear}.  Using the phrase `with constants  $\kappa_{F} > 0$ and $\kappa_{G} > 0$' is accomplished by setting $\kappa_F = \max \{\kappa_{F_i} : i=1, 2, ..., m\}$ and $\kappa_G = \max \{\kappa_{G_i} : i=1, 2, ..., m\}$.

Other notation in this paper will follow that of \cite{rockafellar-wets-1998}.  In particular, for a function $f : \R^n \rightarrow \R$ at a point $\bar{x} \in \dom(f)$, we define the {\em (regular) subdifferential} $\partial f(\bar{x})$ by
    $$\partial f(\bar{x}) = \{ v : f(x) \geq f(\bar{x}) + v^\top (x-\bar{x}) + o(\|x-\bar{x}\|), ~\mbox{for all}~x \in \dom(g)\},$$
where $o(\cdot)$ is the `little-oh' function (see \cite[\S 8]{rockafellar-wets-1998}).  The subdifferential represents a set of generalized gradients and plays an important role in optimization.  For example, the subdifferential can be used to check first order optimality \cite[Thm 10.1]{rockafellar-wets-1998}, and determine descent directions \cite[Ex 8.23]{rockafellar-wets-1998}.

Of interest to this paper is the {\em chain rule} for computing subdifferentials.  In particular, if $f$ is the composition of a convex lsc function $g$ and a smooth vector-valued function $F$, $f = g \circ F$, and $F(\bar{x}) \in \mathrm{int}(\dom(g))$, then the chain rule applies at $\bar{x}$,
    \begin{equation}\label{eq:chainrule}
    \partial f(\bar{x}) = \nabla F(\bar{x})^\top \partial g(F(\bar{x})),
    \end{equation}
see \cite[Thm 10.6]{rockafellar-wets-1998}.

\section{Function Value Approximations}\label{functionvalueapproximations}

Our first result examines the error bounds on the function values of $f = g \circ F$, where $g$ is convex lsc function, $F$ is a smooth vector-valued function, and $F$ is approximated via fully linear models.

\begin{theorem}[Function Value Approximations]\label{thm:fvalues} Suppose
    $$f :\begin{array}{lcl}
    \R^n &\rightarrow& \R \cup \{\infty\} \\
    f(x) &\mapsto& g(F (x)),
    \end{array}$$
where $g : \R^m \rightarrow \R  \cup \{\infty\}$ is convex lsc and $F : \R^n \rightarrow \R^m$ is a smooth vector-valued function.  Suppose $\bar{x} \in \R^n$  with $F(\bar{x}) \in \mathrm{int}(\dom(g))$.  Suppose $\bar{\Delta}>0$ and $\{\tilde F_\Delta : \Delta \in (0, \bar{\Delta})\}$ are fully linear models of $F$ at $\bar{x}$ with constants $\kappa_{F} > 0$ and $\kappa_{G} > 0$.  Define $\tilde{f}_\Delta = g \circ \tilde F_\Delta$.  Then, there exists $\kappa_f$ and $\bar{\Delta}_f > 0$ such that
	\begin{equation}\label{eq:functionerror}
		|f(y) - \tilde{f}_\Delta(y)| \leq \kappa_f \, \Delta^2 \qquad \mbox{ for all }~\Delta \in (0, \bar{\Delta}_f), y \in B_\Delta(\bar x).
    \end{equation}
\end{theorem}

\Proof  As $F(\bar{x}) \in \mathrm{int}(\dom(g))$ and $F$ is continuous, there exists $\delta > 0$ such that $F(y) \in \mathrm{int}(\dom(g))$ for all $y \in \cl(B_{\delta}(\bar{x}))$. Define $\bar{\Delta}_f = \min\{\delta, \bar{\Delta}\}$.

By the continuity of $F$, $\{z \in \R^m : z =F(y), y \in \cl(B_{\bar{\Delta}_f}(\bar{x}))\} \subseteq \mathrm{int}(\dom(g))$ is compact.  As such, $g$ is Lipschitz continuous relative to this set \cite[Cor 8.32]{BauschkeCombettes2011}.  Let $L$ be the Lipschitz constant of $g$ relative to $\{z \in \R^m : z =F(y), y \in \cl(B_{\bar{\Delta}_f}(\bar{x}))\}$.

Given any $\Delta \in (0, \bar{\Delta}_f)$ and $y \in B_\Delta(\bar x)$, we find that
    $$\begin{array}{rcll}
    |f(y) - \tilde{f}_\Delta(y)| &=& |g(F(y)) - g(\tilde F_\Delta(y))| \\
    &\leq& L \|F(y) - \tilde F_\Delta(y)\|\\
    &\leq& L \sum_{i=1}^m |F_i(y) - \tilde F_{i,\Delta}(y)|  \\
    &\leq& L \sum_{i=1}^m \kappa_{F} \Delta^2 \\
    &\leq& L m \kappa_F \Delta^2.
    \end{array}$$
Thus $\kappa_f = L m \kappa_F$ and $\bar{\Delta}_f$ defined above, provide the desired scalars. \qed\medskip

The scalar $\kappa_f$ in Theorem \ref{thm:fvalues} depends on the locally Lipschitz constant of $g$, and on the dimension $m$.  It is not difficult to create examples showing that these values are required to properly define $\kappa_f$.  (Indeed, if $g(z) = \sum_{i=1}^m L |z_i|$ and $|F_i(y) - \tilde F_{i,\Delta}(y)| = \kappa_{F_i} \Delta^2$, then the proof becomes tight.)  However, in the case where $g$ is the max-function, the local Lipschitz constant reduces to $L=1$, and it is possible to remove the dependence on the dimension $m$.  In doing this, a cleaner bound is achieved.

\begin{proposition} Suppose
    $$f :\begin{array}{lcl}
    \R^n &\rightarrow& \R \\
    f(x) &\mapsto& \max(F(x)),
    \end{array}$$
where $F : \R^n \rightarrow \R^m$ is a smooth vector-valued function.  Suppose $\bar{x} \in \R^n$. Suppose $\bar{\Delta}>0$ and $\{\tilde F_\Delta : \Delta \in (0, \bar{\Delta})\}$ are fully linear models of $F$ at $\bar{x}$ with constants $\kappa_{F} > 0$ and $\kappa_{G} > 0$.  Then
	\begin{equation}\label{eq:functionerrormax}
	|f(y) - \tilde{f}_\Delta(y)| \leq \kappa_F \, \Delta^2 \qquad \mbox{ for all }~\Delta \in (0, \bar{\Delta}), y \in B_\Delta(\bar x).
    \end{equation}
\end{proposition}

\Proof Select any $\Delta \in (0, \bar{\Delta})$ and $y \in B_\Delta(\bar x)$.  Let $I \in \argmax \{F_1(y), F_2(y), ..., F_m(y)\}$ and $\tilde{I} \in \argmax \{\tilde{F}_{1,\Delta}(y), \tilde{F}_{2,\Delta}(y), ..., \tilde{F}_{m,\Delta}(y)\}$.

Notice that
    $$\begin{array}{rcll}
    f(y) - \tilde{f}_\Delta(y) &=& \max \{F_1(y), F_2(y), ..., F_m(y)\} - \tilde{F}_{\tilde{I},\Delta}(y) \\
    &\geq& F_{\tilde{I}}(y) - \tilde{F}_{\tilde{I},\Delta}(y)\\
    &\geq&  - \kappa_F \Delta^2,
    \end{array}$$
where the last inequality results from $\tilde{F}_{\tilde{I},\Delta}$ being fully linear models of $F_{\tilde{I}}$.  Similarly,
$$\begin{array}{rcll}
    \tilde{f}_\Delta(y) - f(y) &=&  \max \{\tilde{F}_{1,\Delta}(y), \tilde{F}_{2,\Delta}(y), ..., \tilde{F}_{m,\Delta}(y)\} - F_{I}(y) \\
    &\geq& \tilde{F}_{I,\Delta}(y) - F_{I}(y)\\
    &\geq&  - \kappa_F \Delta^2.
\end{array}$$
Combined, these yield equation \eqref{eq:functionerrormax}.\qed\medskip

We end this section with an example demonstrating the importance of the assumption $F(\bar{x}) \in \mathrm{int}(\dom(g))$ in Theorem \ref{thm:fvalues}.

\begin{example}[Importance of $F(\bar{x}) \in \mathrm{int}(\dom(g))$] \label{ex:importanceofint1} Consider
    $$f :\begin{array}{lcl}
    \R^2 &\rightarrow& \R \\
    f(x) &\mapsto& g(F(x)),
    \end{array}$$
where
    $$g(z) = \left\{\begin{array}{rl}
    0 & \mbox{if } z_1 \geq 0 \\
    \infty & \mbox{otherwise}
    \end{array}\right.$$
and $F(x) = (x_1, x_2)$.  Clearly $g$ is convex lsc and $F$ is smooth.  Let $\tilde{F}_{\Delta}(x) = (x_1 - \Delta^2, x_2)$, and notice that $\tilde{F}_{\Delta}$ are fully linear models of $F$ at $(0, 0)$ with constants $\kappa_F=1$, $\kappa_G = 1$, and $\bar{\Delta}=1$.  However, $F(0,0) \notin \mathrm{int}(\dom(g))$.

Define $\tilde{f}_\Delta(x) = g(\tilde{F}_{\Delta}(x))$, and notice that
    $$|f(0,0) - \tilde{f}_\Delta(0,0)| = |g(0,0) - g(-\Delta^2, 0)| = \infty.$$
Hence, the bound in equation \eqref{eq:functionerror} can never be achieved.\qed

\end{example}

\section{Subdifferential Approximations}\label{subdifferentialapproximations}

We now turn our attention to the error bounds on the subgradient vectors of $f = g \circ F$, where $g$ is convex lsc function, $F$ is a smooth vector-valued function, and $F$ is approximated via fully linear models. Notice that Theorem \ref{thm:gvalues} makes one additional assumption, $F(\bar{x}) = \tilde{F}_\Delta(\bar{x})$, and the results of Theorem \ref{thm:gvalues} are focused only at the point $\bar{x}$ (instead of all points within $\Delta$ of $\bar{x}$.  The need for these assumptions is explored in Examples \ref{ex:importanceFequaltildeF} and \ref{ex:importantofbeingatxbar}.

\begin{theorem}[Subdifferential Approximations]\label{thm:gvalues} Suppose
    $$f :\begin{array}{lcl}
    \R^n &\rightarrow& \R \cup \{\infty\} \\
    f(x) &\mapsto& g(F (x)),
    \end{array}$$
where $g : \R^m \rightarrow \R \cup \{\infty\}$ is convex lsc and $F : \R^n \rightarrow \R^m$ is a smooth vector-valued function.  Suppose $\bar{x} \in \R^n$ and $F(\bar{x}) \in \mathrm{int}(\dom(g))$.  Suppose $\bar{\Delta}>0$ and $\{\tilde F_\Delta : \Delta \in (0, \bar{\Delta})\}$ are fully linear models of $F$ at $\bar{x}$ with constants $\kappa_{F} > 0$ and $\kappa_{G} > 0$.   Suppose $F(\bar{x}) = \tilde{F}_\Delta(\bar{x})$.  Define $\tilde{f}_\Delta = g \circ \tilde F_\Delta$.  Then, there exists $\kappa_g$ such that for all $\Delta \in (0, \bar{\Delta})$:
    \begin{enumerate}
        \item given any $v \in \partial f(\bar{x})$ there exists $\tilde{v} \in \partial \tilde{f}_\Delta(\bar{x})$ such that $\|v - \tilde{v}\| \leq \kappa_g \Delta$, and
        \item given any  $\tilde{v} \in \partial \tilde{f}_\Delta(\bar{x})$ there exists $v \in \partial f(\bar{x})$  such that $\|v - \tilde{v}\| \leq \kappa_g \Delta$.
    \end{enumerate}
\end{theorem}

\Proof Since $g$ is convex lsc, $F$ is smooth, and $F(\bar{x}) \in \mathrm{int}(\dom(g))$, by \cite[Thm 10.6]{rockafellar-wets-1998} we have that
    $$\partial f(\bar{x}) = \nabla F(\bar{x})^\top \partial g(F(\bar{x})).$$
Similarly,
    $$\partial \tilde{f}_\Delta(\bar{x}) = \nabla \tilde{F}_\Delta(\bar{x})^\top \partial g(\tilde{F}_\Delta(\bar{x})).$$
Applying $F(\bar{x}) = \tilde{F}_\Delta(\bar{x})$, we have that
    $$\partial \tilde{f}_\Delta(\bar{x}) = \nabla \tilde{F}_\Delta(\bar{x})^\top \partial g(F(\bar{x})).$$

Since $F(\bar{x}) \in \mathrm{int}(\dom(g))$, the set $\partial g(F(\bar{x}))$ is bounded \cite[Prop 16.14]{BauschkeCombettes2011}. 
Let $M = \sup \{ \|w\| : w \in  \partial g(F(\bar{x})) \}$.

Selecting any $v \in \partial f(\bar{x})$, there exist $w \in \partial g(F(\bar{x}))$ such that $v = \nabla F(\bar{x})^\top w$.  Define $\tilde{v} = \nabla \tilde{F}_\Delta(\bar{x})^\top w \in \partial \tilde{f}_\Delta(\bar{x})$, and notice that
    \begin{equation}\label{eq:subgradinequalities}
    \begin{array}{rcl}
    \|v - \tilde{v}\| &=& \|\nabla F(\bar{x})^\top w - \nabla \tilde{F}_\Delta(\bar{x})^\top w\|\\
    &=& \left\| \left[\begin{array}{c}\nabla F_1(\bar{x})^\top w - \nabla \tilde{F}_{1,\Delta}(\bar{x})^\top w \\
    \nabla F_2(\bar{x})^\top w - \nabla \tilde{F}_{2,\Delta}(\bar{x})^\top w\\
    \vdots \\
    \nabla F_m(\bar{x})^\top w - \nabla \tilde{F}_{m,\Delta}(\bar{x})^\top w
    \end{array}\right]\right\|\\
    &=&\left( \sum_{i=1}^m (F_m(\bar{x})^\top w - \nabla \tilde{F}_{m,\Delta}(\bar{x})^\top w)^2  \right)^{\frac{1}{2}} \\
    &\leq& \left( \sum_{i=1}^m \|F_m(\bar{x}) - \nabla \tilde{F}_{m,\Delta}(\bar{x})\|^2 \|w\|^2  \right)^{\frac{1}{2}} \\
    &\leq& \left( \sum_{i=1}^m (\kappa_G \Delta)^2 M ^2  \right)^{\frac{1}{2}} \\
    &\leq& M \sqrt{m} \kappa_G \, \Delta.
    \end{array}
    \end{equation}
Conversely, selecting any $\tilde{v} \in \partial \tilde{f}_\Delta(\bar{x})$, there exist $w \in \partial g(F(\bar{x}))$ such that $v = \nabla \tilde{F}_\Delta(\bar{x})^\top w$.  Define $v = \nabla F(\bar{x})^\top w \in \partial f(\bar{x})$, and notice that the sequence of inequalities \eqref{eq:subgradinequalities} again holds.  Therefore, $\kappa_g = M \sqrt{m} \kappa_G$ provides the desired scalar. \qed\medskip

Similar to Theorem \ref{thm:fvalues}, the scalar $\kappa_g$ in Theorem \ref{thm:fvalues} depends on the dimension $m$, and on a constant {\em related} to the Lipschitz constant of $g$.  Indeed, the {\em Lipschitz modulus} of $g$ at $\bar{z}$ is defined
    $$\mathrm{lip} g(\bar{z}) = \lim_{\Delta \searrow 0} \sup_{\stackrel{z', z \in B_\Delta(\bar{z}),}{z' \neq z}} \frac{|f(z') - f(z)|}{\|z' - z\|}.$$
When $g$ is convex and $\bar{z} \in \mathrm{int}(\dom(g))$, then the Lipschitz modulus can be found via
    $$\mathrm{lip} g(\bar{z}) = \sup \{ \|w\| : w \in  \partial g(\bar{z}) \},$$
see \cite[Thm 9.13]{rockafellar-wets-1998}, which is exactly the constant used in the proof of Theorem \ref{thm:gvalues} (with $\bar{z} = F(\bar{x})$).

In the case where $g$ is the max-function, the maximum size of a subgradient vector is $M=1$ and it is possible to remove the dependence on the dimension $m$.  This can be found in \cite[Lem 2.1]{HareNutini2013}, so we do not repeat the result here.

In the case where $g$ is the $\ell_1$ norm, the maximum size of a subgradient vector is $M=1$ and the results reconstruct \cite[Lem 2]{LarsonMenickellyWild2015} (although it should be noted that their notation differs slightly from the notation herein).

The importance of $F(\bar{x}) \in \mathrm{int}(\dom(g))$ to Theorem \ref{thm:gvalues} is clear from the fact that if $F(\bar{x}) \in \mathrm{bndry}(\dom(g))$, then $\partial g(F(\bar{x}))$ is either empty or unbounded \cite[Prop 16.14]{BauschkeCombettes2011}.  If $\partial g(F(\bar{x}))$ is empty, then the theorem becomes moot; while if $\partial g(F(\bar{x}))$ is unbounded, the next example shows what can go wrong.

\begin{example}[Importance of $F(\bar{x}) \in \mathrm{int}(\dom(g))$]\label{ex:importanceofint2} Consider $f$, $g$, and $F$ as defined in Example \ref{ex:importanceofint1}.  At $\bar{x}=(0,0)$, we have $f(0,0) = 0$, $\partial f(0,0) = \{ v : v_1 \leq 0, v_2=0\}$, and $F(0,0) \notin \mathrm{int}(\dom(g))$.

Let $\tilde{F}_{\Delta} = (x_1 - \Delta x_2, x_2)$, and notice that $\tilde{F}_{\Delta}$ are fully linear models of $F$ at $(0, 0)$ with constants $\kappa_F=1$, $\kappa_G = 1$, and $\bar{\Delta}=1$.  Define $\tilde{f}_\Delta(x) = g(\tilde{F}_{\Delta}(x))$, and notice that
    $$\tilde{f}_\Delta(x) = \left\{ \begin{array}{rl}
    0 & \mbox{if}~x_2 \geq \frac{1}{\Delta}x_1 \\
    \infty & \mbox{otherwise}
    \end{array}\right.$$
As such, $\partial \tilde{f}_\Delta(0,0) = \{ v : v_1 \leq 0, v_2=-\Delta v_1 \}$.  In particular, given any $\Delta >0$, the vector $\tilde{v}_{\Delta} = (-1/\Delta, 1)^\top \in \partial \tilde{f}_\Delta(0,0)$.  Moreover, this vector has
    $$\mathrm{dist} (\tilde{v}_{\Delta}, \partial f(0,0)) = 1,$$
which does not converge to $0$ as $\Delta \rightarrow 0$.\qed

\end{example}

Theorem \ref{thm:gvalues} also introduced the new condition that $F(\bar{x}) = \tilde{F}_\Delta(\bar{x})$.  Most common methods for constructing fully linear models satisfy this assumption: for example, linear interpolation \cite[\S 2.3]{conn-scheinberg-vicente-2009}, quadratic interpolation \cite[\S 3.4]{conn-scheinberg-vicente-2009}, and minimum Frobenius norm models \cite[\S 5.3]{conn-scheinberg-vicente-2009}, all satisfy this assumption.  However, note that models based on linear regression \cite[\S 2.3]{conn-scheinberg-vicente-2009} may not satisfy this assumption.  The next example explores what can go wrong when this assumption is removed.

\begin{example}[Importance of $F(\bar{x}) = \tilde{F}_\Delta(\bar{x})$]\label{ex:importanceFequaltildeF} Consider
    $$f :\begin{array}{lcl}
    \R^2 &\rightarrow& \R \\
    f(x) &\mapsto& g(F(x)),
    \end{array}$$
where
    $$g(z) = |x_1| + |x_2|.$$
and $F(x) = (x_1, x_2)$.  Clearly $g$ is convex lsc, $F$ is smooth, and $\dom(g) = \R^2$ so $F(x) \in \mathrm{int}(\dom(g))$ for any $x$.  At $\bar{x}=(0,0)$, we have $f(0,0) = 0$, $\partial f(0,0) = \{ v : -1 \leq v_1 \leq 1, -1 \leq v_2 \leq 1\}$.

Let $\tilde{F}_{\Delta} = (x_1 + \Delta^2, x_2 + \Delta^2)$, and notice that $\tilde{F}_{\Delta}$ are fully linear models of $F$ at $(0, 0)$ with constants $\kappa_F=1$, $\kappa_G = 1$, and $\bar{\Delta}=1$.  However, $\tilde{F}_{\Delta}(0,0) \neq F(0,0)$.

Define $\tilde{f}_\Delta(x) = g(\tilde{F}_{\Delta}(x))$, and notice that
    $$\tilde{f}_\Delta(x) = |x_1 + \Delta^2| + |x_2 + \Delta^2|.$$
In particular, $\tilde{f}_\Delta(x)$ is actually differentiable at $\bar{x} = (0,0)$ with $\nabla \tilde{f}_\Delta(0,0) = (1, 1)^\top$.  Hence, for $(-1,-1)^\top \in \partial f(0,0)$,
    $$\mathrm{dist} ((-1,-1), \partial \tilde{f}_\Delta(0,0)) = \sqrt{8},$$
which does not converge to $0$ as $\Delta \rightarrow 0$.\qed

\end{example}

Finally, we point out that Theorem \ref{thm:gvalues} restricts its analysis to approximate subgradients at the point $\bar{x}$. If $y \in B_\Delta(\bar{x})\setminus\{\bar{x}\}$, then the error bounds no longer apply.

\begin{example}[Failure of Theorem \ref{thm:gvalues} away from $\bar{x}$] \label{ex:importantofbeingatxbar} Consider $f$, $g$, and $F$ as in Example \ref{ex:importanceFequaltildeF}.  Let $\tilde{F}_{\Delta} = (x_1 + \Delta^2 x_2, x_2 + \Delta^2 x_1)$, and notice that $\tilde{F}_{\Delta}$ are fully linear models of $F$ at $(0, 0)$ with constants $\kappa_F=1$, $\kappa_G = 1$, and $\bar{\Delta}=1$.  Note  $\tilde{F}_{\Delta}(0,0) = F(0,0) \in \mathrm{int}(\dom(g))$.

Define $\tilde{f}_\Delta(x) = g(\tilde{F}_{\Delta}(x))$, and notice that
    $$\tilde{f}_\Delta(x) = |x_1 + \Delta^2 x_2| + |x_2 + \Delta^2 x_1|.$$
Notice that $\tilde{f}_\Delta$ is differentiable at the point $\hat{x} = (\epsilon, 0)$ for any $\epsilon > 0$, with $\nabla \tilde{f}_\Delta(\epsilon, 0) = (1 + \Delta^2, 1 + \Delta^2)^\top$. However, $f$ is not differentiable at $(\epsilon, 0)$, and instead has $\partial f(\epsilon, 0) = \{ v : v_1 = 1, v_2 \in [-1, 1]\}$.
Hence, for $(1,-1)^\top \in \partial f(\epsilon,0)$,
    $$\mathrm{dist} ((1,-1), \partial \tilde{f}_\Delta(0,0)) = \sqrt{(1+\Delta^2 - 1)^2 + (1+\Delta^2 + 1)^2} = \sqrt{2 + 2 \Delta^2},$$
which does not converge to $0$ as $\Delta \rightarrow 0$.\qed

\end{example}

\section{Conclusions and Extensions}\label{conclusions}

This paper has derived error bounds resulting from the composition of a convex lsc function with a smooth vector-valued function, under the assumption that fully linear models are available for each component of the vector-valued function. If the convex lsc function is the max-function, then the results recreate results from \cite{HareNutini2013}.  While, if the convex lsc function is the $\ell_1$-norm, then the results recreate results from \cite{LarsonMenickellyWild2015}.

In Theorem \ref{thm:fvalues}, we see that the error bound in the function value retains the same order, over the same neighbourhood, as the original error bound.  In Theorem \ref{thm:gvalues}, we see that the error bound for the first order information retains the same order as the original error bound, but only at the single focal point for the fully linear models.  In examining these Theorems and proofs, it should be immediately clear that if {\em fully linear models} are replaced by {\em fully quadratic models}  \cite{conn-scheinberg-vicente-2008} \cite{WildShoemaker2011}, then the results are trivially adaptable.  In particular, if the error bounds for the fully linear models of the vector-valued function are replaced with
    $$\begin{array}{rcll}
		|F_i(y) - \tilde{F}_{i,\Delta}(y)| &\leq & \kappa_F \Delta^{\Gamma} & \qquad \mbox{ for all }~\Delta  \in (0, \bar{\Delta}), y \in B_\Delta(\bar x), \\
		\mbox{ and } \quad
		\|\nabla F_i(y) - \nabla \tilde{F}_{i,\Delta}(y)\| &\leq&  \kappa_G \Delta^{\Upsilon} & \qquad \mbox{ for all }~\Delta  \in (0, \bar{\Delta}), y \in B_\Delta(\bar x),
	\end{array}$$
then the resulting error bounds in Theorems \ref{thm:fvalues} and \ref{thm:gvalues} will become
	$$|f(y) - \tilde{f}_\Delta(y)| \leq L m \kappa_F \, \Delta^{\Gamma} \qquad \mbox{ for all }~\Delta \in (0, \bar{\Delta}_f), y \in B_\Delta(\bar x),$$
and
     $$\mbox{given any}~v \in \partial f(\bar{x})~\mbox{there exists}~\tilde{v} \in \partial \tilde{f}_\Delta(\bar{x}) ~\mbox{such that}~\|v - \tilde{v}\| \leq M \sqrt{m} \kappa_G \Delta^{\Upsilon},$$
     $$\mbox{given any}~\tilde{v} \in \partial \tilde{f}_\Delta(\bar{x})~\mbox{there exists}~v \in \partial f(\bar{x})~\mbox{such that}~\|v - \tilde{v}\| \leq M \sqrt{m} \kappa_G \Delta^{\Upsilon}.$$

\subsection*{{\bf Acknowledgements}}

This research was partially funded by the Natural Sciences and Engineering Research Council (NSERC) of Canada, Discover Grant \#355571-2013.

\end{document}